\documentclass[a4paper,11pt]{article} \usepackage[dvips]{graphicx}
\usepackage{amssymb} 
\usepackage{amsmath} 
\usepackage{hyperref}
 
\newcommand{\bs}{\boldsymbol} \makeatletter 
\newcommand{\Zset}{\mathbb{Z}}
\def\dddots{\mathinner{\mkern1mu\raise\p@
    \hbox{.}\mkern2mu\raise4\p@\hbox{.}\mkern2mu
    \raise7\p@\vbox{\kern7\p@\hbox{.}}\mkern1mu}}
\newtheorem{theorem}{Theorem} 
\newtheorem{Lemma}{Lemma}
\newtheorem{Remark}{Remark} 
\newtheorem{Corollary}{Corollary}
\newtheorem{Proposition}{Proposition} 
 
\newtheorem{Conjecture}{Conjecture}
\newenvironment{AMS}{\small\bf 2010 AMS subject classification: }{}

 \newfont{\BBB}{msbm10 scaled \magstep1} 
\newfont{\BBS}{msbm10}

\begin{document} 
\title{\bf Trivariate polynomial approximation\\on Lissajous curves   
$\,$\thanks{Supported  
the ex-60$\%$ funds and by the biennial project CPDA124755 
of the University of 
Padova, and by the INdAM GNCS.}}

\author{{\bf L. Bos $^{1}$, S. De Marchi $^{2}$ and M. Vianello 
$^{2}$}}   

\maketitle

\footnotetext[1]{Department od Computer Science, University of Verona 
(Italy).}

\footnotetext[2]{Department of Mathematics, University of Padova (Italy).}

\begin{abstract} 
We study Lissajous curves in the 3-cube, that generate algebraic cubature 
formulas on a special family of rank-1 Chebyshev lattices. These formulas 
are used to construct trivariate hyperinterpolation 
polynomials via a single 1-d Fast Chebyshev Transform (by the 
Chebfun package), and to 
compute discrete extremal sets of Fekete and Leja type for trivariate 
polynomial interpolation. 
Applications could arise in the framework of Lissajous sampling for 
MPI (Magnetic Particle Imaging).
\end{abstract}

\vskip0.5cm \noindent \begin{AMS} {\rm 41A05, 41A10, 41A63, 65D05.} 
\end{AMS} \vskip0.2cm \noindent {\small{\bf Keywords:} 
three-dimensional Lissajous curves, Lissajous sampling, Chebyshev 
lattices, trivariate polynomial interpolation and hyperinterpolation, 
discrete extremal sets.}

\section{Introduction}
During the last decade, a new family of points for bivariate 
polynomial interpolation has been proposed and extensively studied, 
namely the so-called ``Padua points'' of the square; cf. 
\cite{BCDMVX06,BDMVX07,CDMSV11,CDMV05}. They are the 
first known optimal nodal set for total-degree multivariate polynomial 
interpolation, with a Lebesgue constant increasing like $(\log n)^2$, 
$n$ being the polynomial degree. 

One of the key features of the Padua points, 
that has allowed to 
construct the interpolation formula, is that they lie on a suitable 
Lissajous curve, such that the integral of any polynomial of degree $2n$ 
along the curve is equal to the 2d-integral with respect to the product 
Chebyshev measure. More specifically, the Padua points are side contacts 
and self-intersections of the Lissajous curve.

Motivated by that construction, in the present paper we try to extend the 
Lissajous curve technique in dimension 3. Since the resulting curve is not 
self-intersecting, we cannot obtain total-degree polynomial interpolation. 
On the other hand, we are able to generate an algebraic cubature formula 
for the product Chebyshev measure, whose nodes lie on the Lissajous curve 
thus forming a rank-1 Chebyshev lattice (on Chebyshev lattices cf., e.g., 
\cite{CP11}). 

By such formula we can perform polynomial hyperinterpolation, 
which is a discretized orthogonal polynomial expansion \cite{S95}, 
and can be constructed by a single 
1-dimensional Fast Chebyshev Transform 
along the curve. 
Moreover, since the underlying Chebyshev lattices turn out to be Weakly 
Admissible Mehes for total-degree polynomials (cf. \cite{BDMSV11}), we can 
extract 
from them suitable discrete extremal sets of Fekete and Leja type 
for polynomial interpolation (cf. \cite{BDMSV10}). We provide a Matlab 
implementation of the
hyperinterpolation and interpolation scheme, and show some numerical 
examples.  
Applications could arise within the emerging field of MPI (Magnetic 
Particle 
Imaging), cf. \cite{KB13}. 

\section{3d Lissajous curves and Chebyshev lattices}

Below, we shall denote the product Chebyshev measure in
$[-1,1]^3$ by
\begin{equation} \label{chebmeas}
d\lambda=w(\bs{x})d\bs{x}\;,\;\;w(\bs{x})=
\frac{1}{\sqrt{(1-x_1^2)(1-x_2^2)(1-x_3^2)}}\;.
\end{equation}
Moreover, $\mathbb{P}_k^3$ will denote the space of trivariate polynomials 
of degree not exceeding $k$, whose dimension is 
$\mbox{dim}(\mathbb{P}_k^3)=(k+1)(k+2)(k+3)/6$. 

Following the lines of the construction of the Padua points, the strategy 
adopted is to seek a Lissajous curve such that the integral of a 
polynomial in $\mathbb{P}_{2n}^3$ with respect to the Chebyshev 
measure $d\lambda$ is equal (up to a constant 
factor) to the integral of the polynomial along the curve. To this 
purpose, the following integer arithmetic result plays a key role.   

\begin{theorem}
Let be $n\in \mathbb{N}^+$ and $(a_n,b_n,c_n)$ be the integer triple
\begin{equation} \label{triples}
(a_n,b_n,c_n)=\left\{
\begin{array}{ll}
\left(\frac{3}{4}n^2+\frac{1}{2}n,\frac{3}{4}n^2+n,
\frac{3}{4}n^2+\frac{3}{2}n+1\right)
\;,\;\;n\;\mbox{{\em even\/}}
\\ \\
\left(\frac{3}{4}n^2+\frac{1}{4},\frac{3}{4}n^2+\frac{3}{2}n-\frac{1}{4},
\frac{3}{4}n^2+\frac{3}{2}n+\frac{3}{4}\right)
\;,\;\;n\;\mbox{{\em odd\/}}
\end{array}
\right.
\end{equation}
Then, for ever integer triple $(i,j,k)$, not all 
$0$, with $i,j,k\geq 0$ and 
$i+j+k\leq m_n=2n$, we have the property that $i a_n\neq j b_n+k c_n$, $j 
b_n\neq i 
a_n+k c_n$, $k 
c_n\neq i a_n +j b_n$. Moreover, $m_n$ is maximal, in the sense that 
there exist a triple $(i^\ast,j^\ast,k^\ast)$, 
$i^\ast+j^\ast+k^\ast=2n+1$, 
that does not satisfy the property.  
\end{theorem}
\vskip0.5cm
\noindent
{\em Proof.\/} See the Appendix.  
\vskip0.5cm

\begin{Proposition}
Consider the Lissajous curves in $[-1,1]^3$ defined by
\begin{equation} \label{lissa}
\bs{\ell}_n(\theta)=(\cos(a_n \theta),\cos(b_n \theta),\cos(b_n 
\theta))
\;,\;\;\theta\in [0,\pi]\;,
\end{equation}
where $(a_n,b_n,c_n)$ is the sequence of integer triples (\ref{triples}). 

Then, for every total-degree polynomial $p\in \mathbb{P}_{2n}^3$ 
\begin{equation} \label{intcurve}
\int_{[-1,1]^3}{p(\bs{x})\,w(\bs{x})d\bs{x}}
=\pi^2\,\int_0^\pi{p(\bs{\ell}_n(\theta))\,d\theta}\;.
\end{equation}
\end{Proposition}
\vskip0.5cm   
\noindent
{\em Proof.\/} It is sufficient to prove the identity for a polynomial 
basis. Take the total-degree product Chebyshev basis 
$T_i(x_1)T_j(x_2)T_k(x_3)$, $i,j,k \geq 0$, $i+j+k\leq 2n$. 
For $i=j=k=0$, (\ref{intcurve}) is clearly true. For $i+j+k>0$, 
by orthogonality of the basis 
$$
\int_{[-1,1]^3}{T_i(x_1)T_j(x_2)T_k(x_3)\,w(\bs{x})d\bs{x}}=0\;.
$$
On the other hand, 
$$
\int_0^\pi{T_i(\cos(a_n \theta))T_j(\cos(b_n \theta))
T_k(\cos(c_n \theta))\,d\theta}
$$
$$
=\int_0^\pi{\cos(ia_n \theta)\,\cos(jb_n \theta)\,\cos(kc_n 
\theta)\,d\theta}
$$
$$
=\frac{1}{4}\left\{\left. \left. \frac{\sin((ia_n-jb_n-kc_n)\theta)}
{ia_n-jb_n-kc_n}\right|_0^\pi
+\frac{\sin((ia_n+jb_n-kc_n)\theta)}{ia_n+jb_n-kc_n}\right|_0^\pi
\right.
$$
$$
\left. \left. \left.
+\frac{\sin((ia_n-jb_n+kc_n)\theta)}{ia_n-jb_n+kc_n}\right|_0^\pi
+\frac{\sin((ia_n+jb_n+kc_n)\theta)}{ia_n+jb_n+kc_n}\right|_0^\pi\right\}
$$

Now, the fourth summand on the right-hand side is zero since 
$ia_n+jb_n+kc_n>0$, and thus the whole right-hand side is zero if (and 
only if) 
$ia_n-jb_n-kc_n\neq 0$, $ia_n+jb_n-kc_n\neq 0$, $ia_n-jb_n+kc_n\neq 0$, 
which is true by Theorem 1 since $i+j+k\leq 2n$. $\;\;\;\square$

\vskip0.5cm

\begin{Corollary}
Let be $p\in \mathbb{P}_{2n}^3$, $\bs{\ell}_n(\theta)$ the Lissajous 
curve 
(\ref{lissa}) and 
\begin{equation} \label{nu}
\nu=n\,\max\{a_n,b_n,c_n\}=\left\{   
\begin{array}{ll}
\frac{3}{4}n^3+\frac{3}{2}n^2+n
\;,\;\;n\;\mbox{{\em even\/}}
\\ \\
\frac{3}{4}n^3+\frac{3}{2}n^2+\frac{3}{4}n
\;,\;\;n\;\mbox{{\em odd\/}}
\end{array}   
\right.
\end{equation}

Then 
\begin{equation} \label{3dquadrature}
\int_{[-1,1]^3}{p(\bs{x})\,w(\bs{x})d\bs{x}}
=\sum_{s=0}^\mu{w_s\,p(\bs{\ell}_n(\theta_s))}\;,
\end{equation}
where 
\begin{equation} \label{w}
w_s=\pi^2\omega_s\;,\;\;s=0,\dots,\mu\;,
\end{equation}
with 
\begin{equation} \label{cheblatt}
\mu=\nu\;,\;\;\theta_s=\frac{(2s+1)\pi}{2\mu+2}\;,
\;\;\omega_s\equiv \frac{\pi}{\mu+1}\;,\;\;s=0,\dots,\mu\;, 
\end{equation}
or alternatively 
$$
\mu=\nu+1\;,\;\;\theta_s=\frac{s\pi}{\mu}\;,\;\;s=0,\dots,\mu\;,
$$
\begin{equation} \label{cheblatt2}
\omega_0=\omega_\mu=\frac{\pi}{2\mu}\;,\;\;
\omega_s\equiv \frac{\pi}{\mu}\;,\;s=1,\dots,\mu-1\;.
\end{equation}

\end{Corollary}
\vskip0.5cm   
\noindent
{\em Proof.\/} 
Observe that by Proposition 1 and the change of variables  
$t=\cos(\theta)$
$$
\int_{[-1,1]^3}{p(\bs{x})\,w(\bs{x})d\bs{x}}
=\pi^2\,\int_0^\pi{p(\bs{\ell}_n(\theta))\,d\theta}
$$
$$
=\pi^2\,\int_{-1}^1{p(T_{a_n}(t),T_{b_n}(t),T_{c_n}(t))\,
\frac{dt}{\sqrt{1-t^2}}}\;,
$$
where $p(T_{a_n}(t),T_{b_n}(t),T_{c_n}(t))$ is a polynomial of degree 
not exceeding 
$$
2\nu=\max\{ia_n+jb_n+kc_n\,,\,i,j,k\geq 0\,,\,i+j+k\leq 
2n\}=2n\,\max\{a_n,b_n,c_n\}\;.
$$

The conclusion follows by using the classical Gauss-Chebyshev 
or Gauss-Chebyshev-Lobatto  
univariate quadrature 
rules, cf. (\ref{cheblatt}) and (\ref{cheblatt2}) respectively, which are 
exact up to 
degree $2\nu+1$ using the $\mu+1$ nodes 
$\tau_s=\cos(\theta_s)$ and weights $\omega_s$, cf., e.g., 
\cite[Ch. 
8]{MH03}.$\;\;\;\square$ 
\vskip0.5cm

\begin{Remark} (Chebyshev lattices).
{\em We observe that $\{\bs{\ell}_n(\theta_s)\}$, $s=0,\dots,\mu$, are  
3-dimensional rank-1 
Chebyshev lattices (for cubature degree of exactness $2n$) in the 
terminology of 
\cite{CP11}. 
As opposite to \cite{CP12}, where Chebsyhev lattices are generated 
heuristically 
by a search algorithm, here we have a formula to generate 
rank-1 Chebyshev lattices for any degree.  
\/}
\end{Remark}

\subsection{Optimal Tuples and Homogeneous Diophantine Equations}

An  algebraic trivariate polynomial of degree $N$ restricted to the 
Lissajou curve
$\bs{\ell}_n(\theta)$ is a trigonometric polynomial of degree 
$N\max\{a_n,b_n,c_n\}=Nc_n.$ Hence it is some interest to have a triple 
for which the conclusion of Theorem 1 holds with its maximum as small as 
possible. Indeed, we conjecture that the triples \eqref{triples} are 
optimal in this sense.

\begin{Conjecture} Suppose that $(a,b,c)$ is a triple of strictly positive 
integers such that $\max\{a,b,c\}<c_n,$ with $c_n$ given by 
\eqref{triples}. Then 
there exists an  integer triple $(i,j,k),$ not all $0,$ and 
$i+j+k\leq 2n$, such that either $i a= j b+k c$, $j 
b= i a+k c$, or $k c =i a +j b$. In other words, the triples 
\eqref{triples} are those satisfying the conclusion of Theorem 1 having 
the minimum maximum.
\end{Conjecture}

We do not as yet have a proof of this conjecture but can provide a lower 
bound for the minimum maximum of such ``good'' triples with the correct 
order of growth in $n.$

First observe that the conditions of the conclusion of Theorem 1 may be 
expressed more succinctly in terms of a homogeneous linear Diophantine 
equation.

\begin{Lemma}
Suppose that $(a,b,c)$ is a triple of strictly positive integers. Then 
there exists an  integer triple $(i,j,k),$ not all $0,$  and 
$i+j+k\leq N$, such that either $i a= j b+k c$, $j 
b= i a+k c$, or $k c =i a +j b$ iff there exists an integer triple 
$(x,y,z)\in\Zset^3$
such that $|x|+|y|+|z|\le N$ and $xa+yb+zc=0.$
\end{Lemma}
\noindent
{\em Proof.\/} If, for example, $i a= j b+k c,$ then $-ia+jb+kc=0$ and we 
may take $x=-i,$ $y=j$ and $z=k.$ On the other hand, if $xa+yb+zc=0$ then 
not all of
$x,$ $y$ and $z$ can have the same sign. There being an odd number of 
them, two of them have the same sign and the other the opposite sign. By 
multiplying by $-1$ if necessary, we assume that the single sign is 
negative. For example, if it is $x$ that is negative, we may write 
$-xa=yb+zc$ and take $i=-x,$ $j=y$ and $k=z.$
$\;\;\;\square$ 

\medskip

The classical Siegel's Lemma (see e.g. \cite[p. 168]{V99}) gives a bound 
on the order of growth of ``small'' solutions of homogeneous linear 
diophantine equations. We may adapt this to our situation to prove

\begin{Lemma} (A version of Siegel's Lemma) Suppose that $1\le 
n\in\Zset_+.$ Suppose further that $a=[a_1,a_2,\ldots,a_d]\in \Zset_+^d$ 
with $a_i>0,$ $1\le i\le d,$ is such that
\[\max\{a\}\le M\]
where
\[ M:=\left\lfloor \frac{1}{n}{n+d \choose d}\right \rfloor -2\quad 
(=O(n^{d-1})).\]
Then there exists $0\neq x\in \Zset^d$ such that $\sum_{i=1}^d|x_i|\le 2n$ 
and
\[ \sum_{i=1}^d x_i a_i=0.\]
\end{Lemma}
\noindent
{\em Proof.\/} Let $S_d\subset \Zset_+^d$ denote the set of {\it positive} 
tuples $0\neq z\in \Zset_+^d$ such that $\sum_{i=1}^dz_i \le n.$
Then  $\#(S_d)={n+d \choose d}-1.$  

Consider the map $F\,:\,\Zset^d\to \Zset$ given by
\[F(z):=\sum_{i=1}^d a_iz_i.\]
Then  $F(S_d)\subset [1,nM]$
and hence 
\[\#(F(S_d))\le nM.\]
But
\[nM=n\left\{ \left\lfloor\frac{1}{n}{n+d \choose n}\right \rfloor 
-2\right\}\le
{n+d \choose d}-2n<{n+d \choose d}-1, \]
i.e., 
\[\#(F(S_d))< \#(S_d).\]
It follows from the Pigeon Hole Principle that there exists two {\it 
different} tuples
$y^{(1)}\neq y^{(2)}\in S_d$ such that
\[F(y^{(1)})=F(y^{(2)}),\]
i.e.,
\[\sum_{i=1}^d a_i(y_i^{(1)}-y_i^{(2)})=0.\]
The tuple $x:=y^{(1)}-y^{(2)}$ has the desired properties.  
$\;\;\;\square$ 

\medskip
In our context it means that the minimum maximum of ``good'' tuples is at 
least 
\[ M:=\left\lfloor \frac{1}{n}{n+d \choose d}\right \rfloor -2\quad 
(=O(n^{d-1})).\]

\section{Hyperinterpolation on Lissajous curves} 

We shall adopt the following notation. We denote the total-degree 
orthonormal basis 
of $P_n^3([-1,1]^3)$ with respect to the Chebyshev product measure 
(\ref{chebmeas}) by 
\begin{equation} \label{cheb3d}
\hat{\phi}_{i,j,k}(\bs{x})=\hat{T_i}(x_1)\hat{T_j}(x_2)\hat{T_k}(x_3)\;,\;i,j,k\geq 
0\;,\;i+j+k\leq n\;,
\end{equation} 
where $\hat{T}_m(\cdot)$ is the normalized Chebyshev polynomial 
of degree $m$
\begin{equation} \label{cheb1d}
\hat{T}_m(\cdot)=
\sigma_m\cos(m\arccos(\cdot))\;,\;\;
\sigma_m=\sqrt{\frac{1+\mbox{sign}(m)}{\pi}}\;,
\;\;m\geq 0\;, 
\end{equation}
with the convention that $\mbox{sign}(0)=0$.

We recall that hyperinterpolation is a discretized expansion of a function 
in series of orthogonal polynomials up to total-degree $n$ on a given 
$d$-dimensional compact region $K$, where the Fourier-like coefficients 
are 
computed  
by a cubature formula exact on $\mathbb{P}^d_{2n}(K)$. 
It was proposed by Sloan in the seminal paper \cite{S95} in order to 
bypass the intrinsic difficulties 
of polynomial interpolation in the multivariate setting, and since then 
has been successfully used in several instances, for example on the sphere 
\cite{HS07}.

Given a function $f\in C([-1,1]^3)$, in view of the algebraic cubature 
formula 
(\ref{3dquadrature}), the hyperinterpolation polynomial of $f$ is 
\begin{equation} \label{hypf}
\mathcal{H}_nf(\bs{x})=\sum_{0\leq i+j+k\leq n}{C_{i,j,k}\,
\hat{\phi}_{i,j,k}(\bs{x})}\;,
\end{equation}
where 
\begin{equation} \label{hypcoeffs}
C_{i,j,k}=\sum_{s=0}^\mu{w_s\,f(\bs{\ell}_n(\theta_s))
\,\hat{\phi}_{i,j,k}(\bs{\ell}_n(\theta_s))}\;.
\end{equation} 

Observe that by construction $\mathcal{H}_nf=f$ for every 
$f\in \mathbb{P}_n^3$, i.e.,  $\mathcal{H}_n$ is a projection operator.  
Among the properties of the hyperinterpolation 
operator, not depending on the specific cubature formula  
provided it is exact up to degree $2n$ for the product Chebyshev measure, 
we recall the following bound for the $L^2$ error, 
\begin{equation} \label{L2err}
\|f-\mathcal{H}_nf\|_2\leq 
2\pi^3\,E_n(f)\;,\;\;E_n(f)=\inf_{p\in 
\mathbb{P}_n}\|f-p\|_\infty\;.
\end{equation}
Consider the uniform operator norm (i.e., the Lebesgue constant)  
\begin{equation} \label{hypnorm}
\|\mathcal{H}_n\|=\sup_{f\neq 0}
\frac{\|\mathcal{H}_nf\|_\infty}{\|f\|_\infty}
=\max_{\bs{x}\in [-1,1]^3}{\sum_{s=0}^\mu{w_s\,
\left|K_n(\bs{x},\bs{\ell}_n(\theta_s))\right|}}\;,
\end{equation}
where $K_n(\bs{x},\bs{y})=\sum_{0\leq i+j+k\leq n}
{\hat{\phi}_{i,j,k}(\bs{x})\hat{\phi}_{i,j,k}(\bs{y})}$ is the reproducing 
kernel of $\mathbb{P}_n^3$ with respect to the product Chebyshev 
measure (\ref{chebmeas}), cf. \cite{DX01}.

In \cite{DMSV14} the bound 
$\|\mathcal{H}_n\|=\mathcal{O}((\sqrt{n})^3)$ 
has been 
obtained,  
as a consequence of a general result connecting multivariate Christoffel 
functions and hyperinterpolation operator norms. On the other hand, 
by proving a conjecture stated in \cite{DMVX09},  
the fine bound  
\begin{equation} \label{WWW}
\|\mathcal{H}_n\|=\mathcal{O}((\log n)^3) 
\end{equation}
has been provided in \cite{WWW14}, 
which corresponds to the minimal growth of a polynomial projection 
operator, in view of \cite{SV09}. Since $\mathcal{H}_n$ is a projection, 
we get the $L^\infty$ error bound
\begin{equation} \label{Linftyerr}
\|f-\mathcal{H}_nf\|_\infty=\mathcal{O}\left((\log n)^3\,E_n(f)\right)\;.
\end{equation}

We show now that the hyperinterpolation 
coefficients $\{C_{i,j,k}\}$ can be computed by a single 1-dimensional 
discrete Chebyshev transform along the Lissajous curve.  

\begin{Proposition}
Let be $f\in C([-1,1]^3)$, $(a_n,b_n,c_n)$ the sequence of integer 
triples (\ref{triples}), and $\nu$, $\mu$, $\{\theta_s\}$, $\omega_s$, 
$\{w_s\}$ 
as 
in Corollary 1. The 
hyperinterpolation coefficients of $f$ generated by (\ref{3dquadrature}) 
can be computed as 
\begin{equation} \label{chebtransf1d}
C_{i,j,k}
=\frac{\pi^2}{4}\,\sigma_{ia_n}\sigma_{jb_n}\sigma_{kc_n}\,
\left(\frac{\gamma_{\alpha_1}}{\sigma_{\alpha_1}}
+\frac{\gamma_{\alpha_2}}{\sigma_{\alpha_2}}
+\frac{\gamma_{\alpha_3}}{\sigma_{\alpha_3}}
+\frac{\gamma_{\alpha_4}}{\sigma_{\alpha_4}}\right)\;,
\end{equation}
$$
\alpha_1=ia_n+jb_n+kc_n\;,\;\;\alpha_2=\left|ia_n+jb_n-kc_n\right|\;,
$$
$$
\alpha_3=\left|ia_n-jb_n\right|+kc_n\;,\;\;\alpha_4=\left|\left|ia_n-jb_n
\right|-kc_n\right|\;,
$$
where $\{\gamma_m\}$ are the first $\nu+1$ coefficients of the 
discretized Chebyshev 
expansion 
of $f(T_{a_n}(t),T_{b_n}(t),T_{c_n}(t))$, $t\in 
[-1,1]$, namely   
\begin{equation} \label{chebcoeffs}
\gamma_m=\sum_{s=0}^\mu{\omega_s\,\hat{T}_m(\tau_s)\,\,
f(T_{a_n}(\tau_s),T_{b_n}(\tau_s),T_{c_n}(\tau_s))}\;,
\end{equation}
$m=0,1,\dots,\nu$, with $\tau_s=\cos(\theta_s)$, $s=0,1,\dots,\mu$.
\end{Proposition}
\vskip0.5cm   
\noindent
{\em Proof.\/} By the change of variables $\theta=\arccos(t)$ which 
gives 
$$
\bs{\ell}_n(\theta)=(T_{a_n}(t),T_{b_n}(t),T_{c_n}(t))\;,
$$ 
and by the 
classical identity 
$T_h(t)T_k(t)=\frac{1}{2}\,\left(T_{h+k}(t)+T_{|h-k|}(t)\right)$ (cf., 
e.g., \cite[\S 2.4.3]{MH03}), we get  
$$
f(\bs{\ell}_n(\theta))
\,\hat{\phi}_{i,j,k}(\bs{\ell}_n(\theta))=f(\bs{\ell}_n(\arccos(t)))
\,\,\hat{T}_{ia_n}(t)\hat{T}_{jb_n}(t)\hat{T}_{kc_n}(t)
$$
$$
=f(\bs{\ell}_n(\arccos(t)))\,\sigma_{ia_n}\sigma_{jb_n}\sigma_{kc_n}
\,\frac{1}{4}\,\left(T_{\alpha_1}(t)
+T_{\alpha_2}(t)+T_{\alpha_3}(t)+T_{\alpha_4}(t)\right)\;,
$$
and hence we have (\ref{chebtransf1d})-(\ref{chebcoeffs})  
in view of (\ref{hypcoeffs}), 
and the fact that $\{\tau_s\}$ are the nodes of the 
univariate Gauss-Chebyshev or Gauss-Chebyshev-Lobatto formula, 
with  
weights 
$\omega_s$, cf. (\ref{cheblatt})-(\ref{cheblatt2}). $\;\;\;\square$

\vskip0.5cm

\begin{Remark} (Lissajous sampling).  
{\em Hyperinterpolation polynomials on $d$-dimensional cubes can be 
constructed by other cubature formulas for the product Chebyshev measure, 
that can be more efficient in terms of number of 
function evaluations required at a given exactness degree. For 
example, a formula of exactness degree $2n$ with $\mathcal{O}(n^4)$ 
nodes for the 3-cube has been provided in \cite{DMVX09}, and used in a 
FFT-based 
implementation of hyperinterpolation. Other formulas, in particular 
Godzina's blending formulas \cite{G95}, that have the 
lowest cardinality known in $d$-dimensional cubes, have been used in 
the package \cite{CP13}. All such formulas are based on Chebyshev lattices 
of rank greater than 1, that are suitable unions of product Chebyshev 
subgrids.  

A first advantage of rank-1 Chebyshev lattices, as observed in 
general in \cite{CP11}, is that a single 1-dimensional FFT is needed 
to compute the hyperinterpolation polynomials. In the present context 
of sampling on Lissajous curves of the 3-cube, this is manifest in 
Proposition 2. 

On the other hand, one of the most most interesting features of 
hperinterpolation on Lissajous curves 
arises in connection with medical imaging applications, in particular 
with the emerging 3d MPI (Magnetic Particle Imaging) technology. 
Indeed, Lissajous sampling 
is one of the most common sampling methods within this technology, 
since it can be generated by suitable electromagnetic fields with 
different frequencies in the components, cf., e.g., \cite{KB13,MWG00}.   
Choosing the frequencies (\ref{triples}) that generate the 
specific 3d Lissajous curves (\ref{lissa}), a clear connection with 
multivariate polynomial approximation comes out, that could be useful 
in the corresponding data processing and analysis.  
\/}
\end{Remark}

\vskip0.5cm
\begin{Remark} (Clenshaw-Curtis type cubature).
{\em
The availability of an hyperinterpolation operator with respect to a given 
density function (here the trivariate Chebyshev density) allows us to 
easily construct algebraic cubature formulas for other densities, 
generalizing the Clenshaw-Curtis quadrature approach (cf., 
e.g., \cite{MH03}). Indeed, if the ``moments''
\begin{equation} \label{moments}
m_{i,j,k}=\int_{[-1,1]^3}{\hat{\phi}_{i,j,k}(\bs{x})\,\xi(\bs{x})d\bs{x}}\;,
\;\;i,j,k\geq 0\;,\;\;i+j+k\leq n 
\end{equation}
are known, where $\xi \in L^1_+((-1,1)^3)$, as shown in \cite{SVZ08} we 
can construct by (\ref{hypcoeffs}) the cubature formula
$$
\int_{[-1,1]^3}{\mathcal{H}_n f(\bs{x})\,\xi(\bs{x})d\bs{x}}
=\sum_{0\leq i+j+k\leq n}{C_{i,j,k}\,m_{i,j,k}}
$$
\begin{equation} \label{CC}
=\sum_{s=0}^\mu{W_s\,f(\bs{\ell}_n(\theta_s))}\;,\;\;
W_s=w_s\,\sum_{0\leq i+j+k\leq 
n}{m_{i,j,k}\,\hat{\phi}_{i,j,k}(\bs{\ell}_n(\theta_s))}\;,
\end{equation}
which is exact for all polynomials in $\mathbb{P}_n^3$. The resulting 
weights $\{W_s\}$ are not all positive, in general, but if $\xi/w\in 
L^2((-1,1)^3)$, which is true for example for the Lebesgue measure 
$\xi(\bs{x})\equiv 1$, it can be proved that 
\begin{equation} \label{stabCC}
\lim_{n\to \infty}{\sum_{s=0}^\mu{|W_s|}}=
\int_{[-1,1]^3}{\frac{\xi(\bs{x})}{w(\bs{x})}\,d\bs{x}}\;, 
\end{equation}
thus ensuring convergence and 
stability of the cubature formula; cf. \cite{SVZ08}. 

We stress that these Clenshaw-Curtis type cubature formulas 
are based on {\em Lissajous sampling} (see Remark 2), and by Proposition 
2 can 
be constructed 
by a single 1-dimensional discrete Chebyshev transform along the Lissajous 
curve (i.e., by a single 1-dimensional FFT).\/}  
\end{Remark}

\vskip0.5cm
\begin{Remark} (Weakly Admissible Meshes and Discrete Extremal Sets).
{\em In the recent literature on multivariate polynomial approximation, 
the notion
of ``Weakly Admissible Mesh'' has emerged as a basic tool, from both the
theoretical and the computational point of view; cf., e.g., 
\cite{BDMSV10,BDMSV11,CL08} and the references therein.

We recall that a Weakly Admissible Mesh (WAM) is a sequence of finite
subsets
of a multidimensional (polynomial-determining) compact set, say
$\mathcal{A}_n\subset K\subset \mathbb{R}^d$ (or $\mathbb{C}^d$), which
are {\em norming
sets}
for total-degree polynomial subspaces,
\begin{equation} \label{wam}
\|p\|_{\infty,K}\leq 
C(\mathcal{A}_n)\,\|p\|_{\infty,\mathcal{A}_n}\;,\;\;\forall p\in
\mathbb{P}_n^d\;,
\end{equation}
where both $C(\mathcal{A}_n)$ and $\mbox{card}(\mathcal{A}_n)$ increase at
most polynomially with
$n$. Here, $\mathbb{P}_n^d$ denotes the space of $d$-variate
polynomials of degree not exceeding $n$, and $\|f\|_{\infty,X}$
the sup-norm of a
function $f$ bounded on the (discrete or continuous) set $X$.
Observe that
necessarily $\mbox{card}(\mathcal{A}_n)\geq
\mbox{dim}(\mathbb{P}_n^d)$.

Among their properties, we quote that
WAMs are preserved by affine transformations,
can be constructed incrementally by finite union and product, and are
``stable" under small perturbations \cite{PV13}. 
It has been shown in the seminal paper \cite{CL08} that WAMs 
are nearly
optimal for polynomial least-squares approximation in
the uniform norm. Moreover, the interpolation Lebesgue constant 
of Fekete-like extremal sets extracted from such 
meshes, say $\mathcal{F}_n$ (that are points maximizing the Vandermonde 
determinant on $\mathcal{A}_n$), has the bound
\begin{equation} \label{lebfek}
\Lambda(\mathcal{F}_n)\leq \mbox{dim}(\mathbb{P}_n^d)\,
C(\mathcal{A}_n)\;.
\end{equation}

Now, the Chebyshev lattices 
\begin{equation} \label{cheblattwam}
\mathcal{A}_n=\{\bs{\ell}_n(\theta_s)\,,\,\,s=0,\dots,\mu\}
\end{equation} 
in (\ref{cheblatt})-(\ref{cheblatt2}), form a WAM for $K=[-1,1]^3$, 
with $C(\mathcal{A}_n)=\mathcal{O}((\log n)^3)$. 
In fact, the corresponding hyperinterpolation operator 
$\mathcal{H}_n$ being a projection on $\mathbb{P}_n^3$, we get 
by (\ref{WWW})  
\begin{equation} \label{cheblattwamest}
\|p\|_{\infty,[-1,1]^3}=\|\mathcal{H}_n p\|_{\infty,[-1,1]^3}
\leq \|\mathcal{H}_n\|\,\|p\|_{\infty,\mathcal{A}_n}=
\mathcal{O}((\log n)^3)\,\|p\|_{\infty,\mathcal{A}_n}\;.
\end{equation} 

In the next Section, we shall use the fact that Fekete-like extremal sets 
extracted from 
$\mathcal{A}_n=\{\bs{\ell}_n(\theta_s)\,,\,\,s=0,\dots,\mu\}$ provide 
a {\em Lissajous sampling} approach to {\em trivariate polynomial 
interpolation}. 
\/}
\end{Remark}
\vskip0.5cm

\section{Implementation and numerical examples} 

\subsection{Hyperinterpolation by Lissajous sampling}
In view of Proposition 2, hyperinterpolation on the Lissajous curve can be 
implemented by a single 1-dimensional Discrete Chebyshev Transform, 
i.e., by a single 1-dimensional FFT. We shall concentrate on sampling 
at the Chebyshev-Lobatto points, since in this case we can 
conveniently resort 
to the powerful {\em Chebfun package\/} (cf. \cite{DHT14}). Sampling at 
the 
Chebyshev zeros can be treated in a similar way.  

Indeed, in view of a well-known discrete orthogonality property 
of the Chebyshev polynomials,  
the interpolation polynomial of a function $g$ 
at the Chebyshev-Lobatto points can be written as 
\begin{equation} \label{interpol}
\pi_\mu(t)=\sum_{m=0}^\mu{c_mT_m(t)}
\end{equation}
where 
$$
c_m=\frac{2}{\mu}\,\sum_{s=0}^\mu{''\,T_m(\tau_s)\,g(\tau_s)}\;,\;\;
m=1,\dots,\mu-1\;,
$$
\begin{equation} \label{chebcoeffs2}
c_m=\frac{1}{\mu}\,\sum_{s=0}^\mu{''\,T_m(\tau_s)\,g(\tau_s)}\;,\;\;m=0,\mu\;,
\end{equation}
the double prime indicating that the first and the last terms of the sum 
have to be halved (cf., e.g., \cite[\S 6.3.2]{MH03}).

Applying this interpolation formula to 
$g(t)=f(T_{a_n}(t),T_{b_n}(t),T_{c_n}(t))$ and comparing with 
the discrete Chebyshev expansion coefficients (\ref{chebcoeffs}), we 
obtain by easy calculations 
\begin{equation} \label{coeffrel}
\frac{\gamma_m}{\sigma_m}=\left\{\begin{array}{ll} 
\frac{\pi}{2}\,c_m\;,\;\;m=1,\dots,\mu-1 \\ \\
\pi\,c_m\;,\;\;m=0,\mu
\end{array}
\right.
\end{equation}
i.e., the 3-dimensional hyperinterpolation 
coefficients (\ref{chebtransf1d}) can be computed 
by the $\{c_m\}$ and (\ref{coeffrel}).

The coefficients of Chebyshev-Lobatto 
interpolation (\ref{chebcoeffs2}) are 
at the core of the Chebfun package, cf. \cite{BT04,T14}. A single call to 
the Chebfun basic 
function 
{\tt chebfun\/} on $f(T_{a_n}(t),T_{b_n}(t),T_{c_n}(t))$, truncated 
at the $(\mu+1)$th-term, produces all the relevant coefficients 
$\{c_m\}$ in an 
extremely fast and stable way. 

For example, by the Matlab code \cite{DMV14} we can compute in about 1 
second the  
$\mu=\frac{3}{4}n^3+\frac{3}{2}n^2+n+2=765102$ coefficients for $n=100$ 
with functions such as 
\begin{equation} \label{testf}
f_1(\bs{x})=\exp(-c\|\bs{x}\|_2^2)\;,\;c>0\;,\;\;
f_2(\bs{x})=\|\bs{x}\|_2^\beta\;,\;\beta>0\;,
\end{equation}  
from which we get by (\ref{chebtransf1d}) the $(n+1)(n+2)(n+3)/6=176851$ 
coefficients of trivariate hyperinterpolation at degree $n=100$.  
All the numerical tests have been made by Chebfun 5.1, in Matlab 7.7.0 
with an Athlon 64 X2 Dual Core 4400+ 2.40GHz processor.

For the purpose of illustration, in Figure 1 we show the relative errors 
(in the Euclidean norm on a suitable control grid) for 
two polynomials of degree 10 and 20, respectively, and for the test  
functions $f_1$ and $f_2$ in (\ref{testf}). Observe the Gaussian $f_1$ is 
analytic, with variation rate determined by the parameter $c$, 
whereas the power function $f_2$ has 
finite regularity, determined by the parameter $\beta$.  

Notice that the error decreases with the degree to a certain threshold
above machine precision and thereafter does not improve. This is likely 
due to the fact that we require the summation of a large number of terms, 
for which a non-negligible error is to be expected. For practical 
applications this is of little import.


In Figures 2 and 3 one can see the Chebyshev lattice on the Lissajous 
curve for polynomial degree $n=5$.

\vskip0.5cm 
\begin{figure}[!ht] \centering
\includegraphics[scale=0.50,clip]{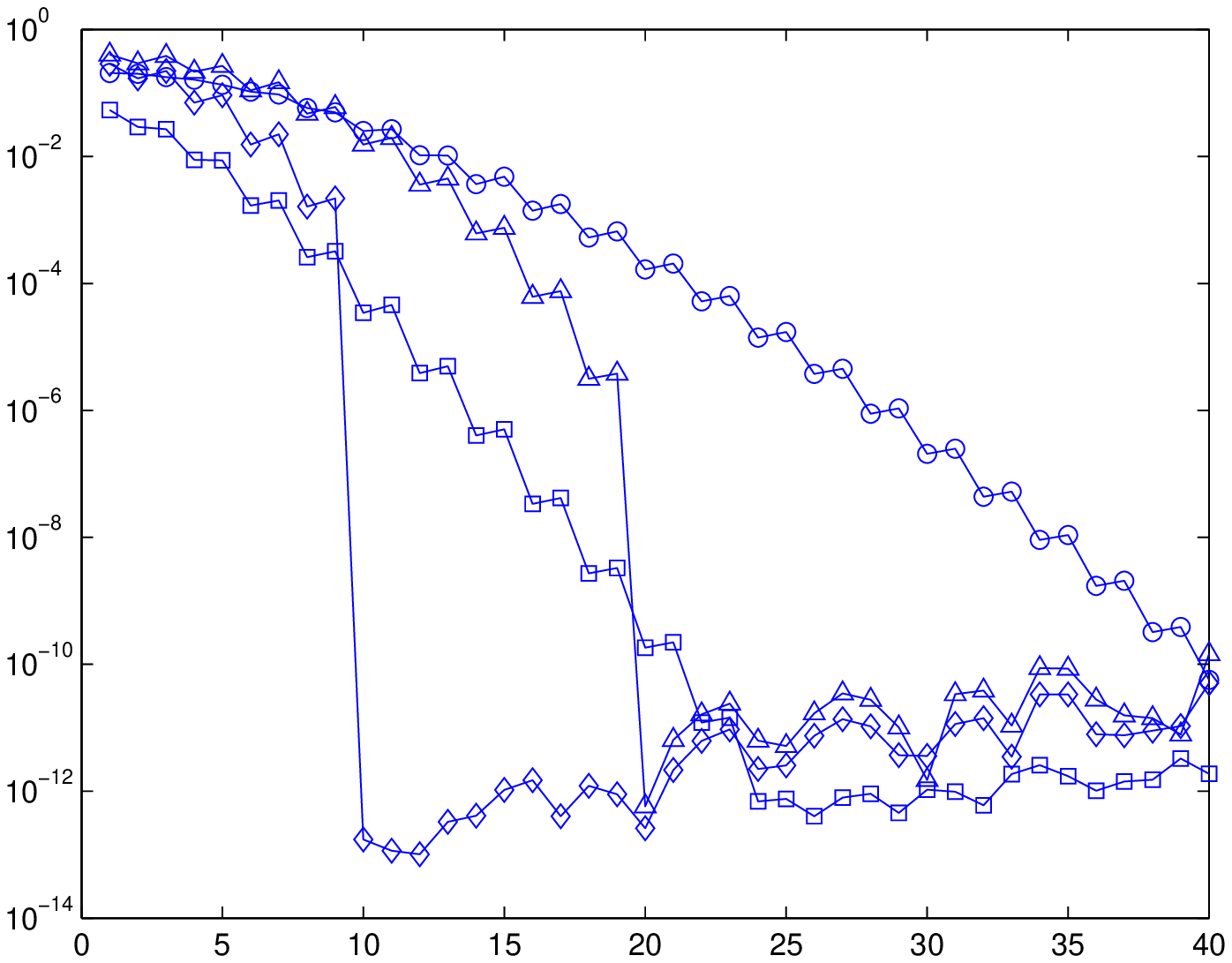}\hfill
\includegraphics[scale=0.50,clip]{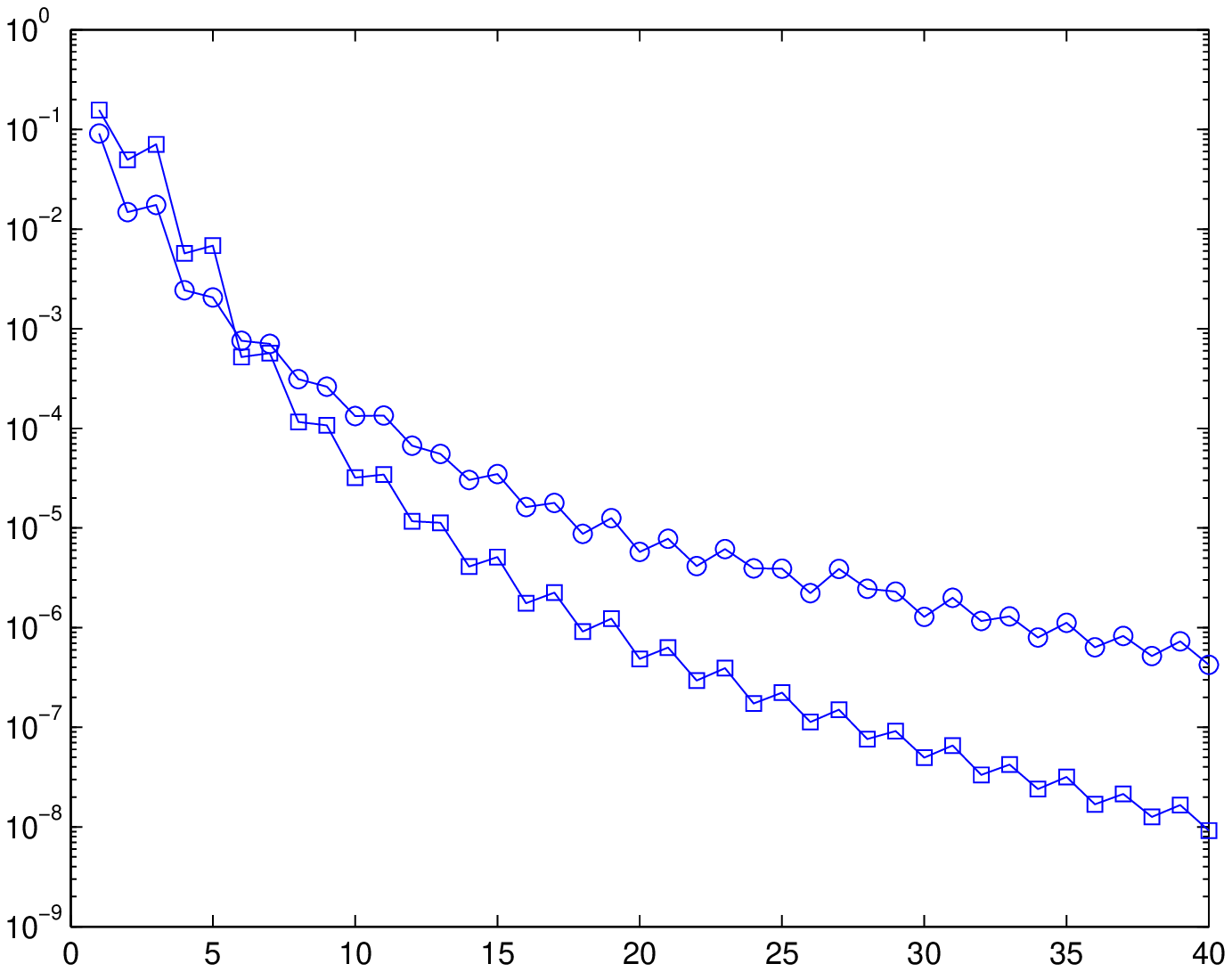}\hfill
\caption{Top: Hyperinterpolation errors for the trivariate polynomials 
$\|\bs{x}\|_2^{2k}$ 
with $k=5$ (diamonds) and $k=10$ (triangles), and for the trivariate 
function $f_1$ 
with $c=1$ (squares) and $c=5$ (circles). Bottom: Hyperinterpolation 
errors for the trivariate function $f_2$ 
with $\beta=5$ (squares) and $\beta=3$ 
(circles).}
\label{hypererr}
\end{figure}
\vskip0.5cm

\subsection{Interpolation by Lissajous sampling} 

Concerning polynomial interpolation in the cube by sampling on the 
Lissajous curve, 
we resort to
the approximate versions of Fekete points (points that maximize
the absolute value of the Vandermonde determinant) 
studied in several recent 
papers \cite{BCLSV11,BDMSV10,SV09}. By (\ref{lebfek}), it makes 
sense to start from a WAM, namely the Chebyshev lattice 
$\mathcal{A}_n$ in  
(\ref{cheblattwam}), by the 
corresponding 
Vandermonde-like matrix 
\begin{equation} \label{chebvand} 
V=V(\mathcal{A}_n;\bs{\phi})\in \mathbb{R}^{M\times N}\;,\;\;
M=\mbox{card}(\mathcal{A}_n)=\mu+1\;,\;\;N=\mbox{dim}(\mathbb{P}_n^3)\;, 
\end{equation}
(cf. (\ref{cheblatt})-(\ref{cheblatt2}) for the definition of $\mu$), 
where 
$$
\bs{\phi}=\{\phi_{i,j,k}\}\;,\;\; 
\phi_{i,j,k}(\bs{x})=T_i(x_1)T_j(x_2)T_k(x_3)\;,\;\;0\leq i+j+k\leq n\;,
$$ 
is the total-degree trivariate Chebyshev 
orthogonal 
basis, suitably ordered (we adopt the {\em graded lexicographical 
ordering}, 
that is the lexicographical ordering within each subset of triples 
$(i,j,k)$ such that $i+j+k=r$, $r=0,\dots,n$). The $(p,q)$ entry of 
$V$ is the $q$-th element of the ordered basis computed in the $p$-th 
element of the nodal array. We recall that the 
choice of the Chebyshev 
orthogonal basis allows to avoid the extreme ill-conditiong of Vandermonde 
matrices in the standard monomial basis.   

The problem of selecting
a $N\times N$ square submatrix with maximal determinant from a given
$M\times N$ rectangular
matrix
is known to be NP-hard \cite{CMI09}, but can be solved in an
approximate way by two simple {\em greedy} algorithms, that are fully
described and analyzed in \cite{BDMSV10}. These algorithms produce
two interpolation nodal sets, called {\em discrete extremal sets}. 

The first, that computes the so-called {\em Approximate Fekete Points}
(AFP),
tries to maximize iteratively submatrix volumes until a maximal volume
$N\times N$ submatrix of $V$ is obtained, and can be based on the famous
{\em QR factorization with column pivoting} \cite{BG65}, applied to   
$V^t$ (that in Matlab is implemented by the matrix left division
or backslash operator, cf. \cite{matlab}).   
See \cite{CMI09}
for the notion of volume generated by a set of vectors, which
generalizes the geometric concept related to parallelograms
and parallelepipeds (the volume and determinant notions coincide
on a square matrix).

The second, that computes the so-called {\em Discrete Leja Points} (DLP),
tries to maximize iteratively submatrix determinants, and is based
simply on {\em Gaussian elimination with row pivoting} applied to 
the Vandermonde-like matrix $V$. 

Denoting by $A$ the $M\times 2$ array of the WAM nodal
coordinates, the
corresponding computational steps, written in a Matlab-like
style,
are
\begin{equation} \label{AFP}
\bs{w}=V\backslash\bs{v};\;\bs{s}={\tt find}(\bs{w}\neq \bs{0});\;
\mathcal{F}_n^{\small{AFP}}=A(\bs{s},:);
\end{equation}
for AFP, where $\bs{v}$ is any nonzero $N$-dimensional vector,
and
\begin{equation} \label{DLP}
[L,U,\bs{\sigma}]={\tt
LU}(V,``{\tt vector}");\;\bs{s}=\bs{\sigma}(1:N);\;
\mathcal{F}_n^{\small{DLP}}=A(\bs{s},:);
\end{equation}
for DLP. In (\ref{DLP}), we refer to the Matlab version of the LU
factorization that produces a row permutation vector. In both
algorithms, we eventually select an index subset
$\bs{s}=(s_1,\dots,s_N)$, that
extracts a Fekete-like discrete extremal set $\mathcal{F}_n$ of the
cube from the WAM
$\mathcal{A}_n$.

Once the underlying extraction WAM   
has been fixed, differently from the continuum Fekete points,
Approximate Fekete Points depend on the
choice of the basis, and Discrete Leja Points depend also on its
order.
An important feature is that Discrete Leja Points form a {\em
sequence}, i.e., if the polynomial basis is such that its first
$N_r=\mbox{dim}(\mathbb{P}^d_r)$ elements span $\mathbb{P}^d_r$,
$1\leq r\leq n$ (as it happens with the graded lexicographical 
ordering of the Chebyshev basis), then the
first $N_r$
Discrete Leja Points are a unisolvent set for interpolation in
$\mathbb{P}^d_r$.

Under the latter assumption for Discrete Leja Points, the two
families of discrete
extremal sets share the same
asymptotic behavior, which by a recent deep result in pluripotential 
theory, cf. \cite{BBWN11}, is exactly that
of the continuum Fekete
points: the corresponding uniform discrete probability measures  
converge weakly to the
{\em pluripotential theoretic equilibrium measure} of the underlying 
compact
set, cf. \cite{BCLSV11,BDMSV10}. In the present case of the cube, 
such a measure is the product Chebyshev measure (\ref{chebmeas}), 
with scaled density $w(\bs{x})/\pi^3$. 

We give now some numerical examples, that can be reproduced by 
the Matlab package \cite{DMPSV14}. First, in  
Figures 2-3 we show the 
Approximate Fekete Points extracted from the Chebyshev lattice on the 
Lissajous curve for degree $n=5$. In Figure 4 we display the numerically 
evaluated Lebesgue constants of the Approximate Fekete Points
and Discrete Leja Points for degree $n=1,2,\dots,30$. For both the 
nodal families, the Lebesgue 
constant turns out to be much lower than the upper bound (\ref{lebfek}), 
and even 
lower than $N=\mbox{dim}(\mathbb{P}_n^3)$, a theoretical upper bound  
for the continuum Fekete points. In particular, the Lebesgue constant of 
Approximate Fekete Points seems to increase quadratically with respect to 
the degree, at least in the given degree range. 

Finally, In Figure 5 we show the relative interpolation errors for the 
two test functions $f_1$ and $f_2$ of Figure 1. Since the Discrete Leja 
Points form a sequence, as discussed above, we have computed them 
once and for all for degree $n=30$, and then used the nested subsequences 
with  
$N_r=\mbox{dim}(\mathbb{P}^d_r)$ elements for interpolation at degree 
$r=1,\dots,30$. The corresponding file of nodal coordinates can be 
downloaded from 
\cite{DMPSV14}. The relevant indexes $(s_1,s_2,\dots,s_{N_{30}})$ 
corresponding to the extraction of  
the Discrete Leja 
Points from the Chebyshev lattice (\ref{cheblattwam})-(\ref{cheblatt2})   
at degree 30, could be used in applications, 
such as MPI 
\cite{KB13}, where a trivariate function is not known or computable 
everywhere, 
but can be sampled just by travelling along the Lissajous curve.

\vskip0.5cm  
\begin{figure}[!ht] \centering
\includegraphics[scale=0.70,clip]{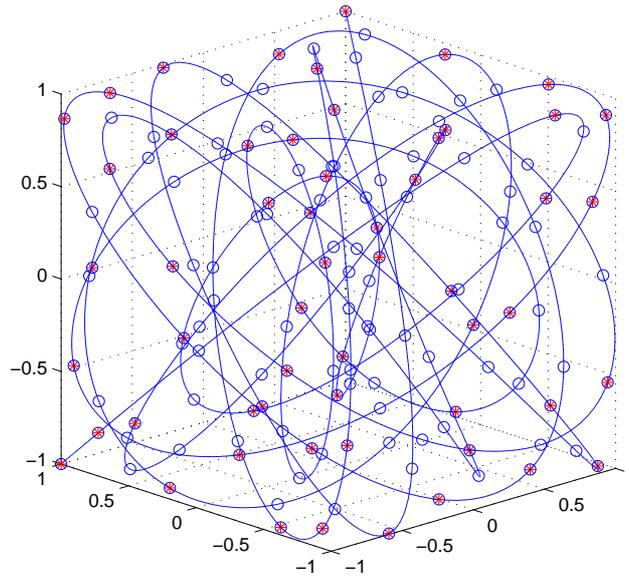}\hfill
\caption{The Chebyshev lattice (circles) and the extracted Approximate 
Fekete Points (asterisks), on the Lissajous curve for polynomial degree 
$n=5$.} \label{feklissa5}
\end{figure}
\vskip0.5cm

\vskip0.5cm
\begin{figure}[!ht] \centering
\includegraphics[scale=0.50,clip]{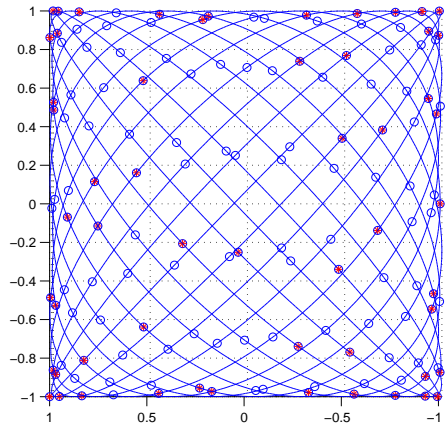}\hfill
\caption{A face projection of the Lissajous curve above with 
the sampling nodes.} 
\label{facelissa5}
\end{figure}
\vskip0.5cm

\vskip0.5cm
\begin{figure}[!ht] \centering
\includegraphics[scale=0.50,clip]{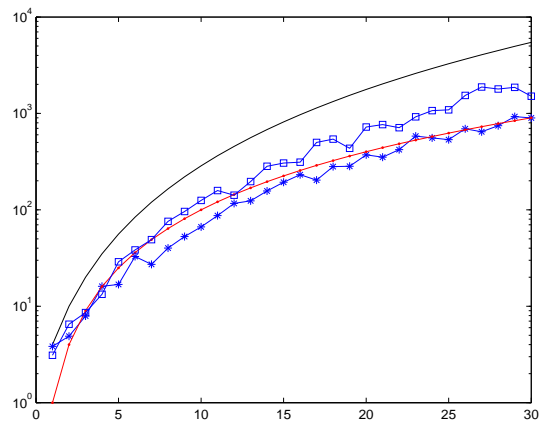}\hfill
\caption{Lebesgue constants (log scale) of the Approximate Fekete Points 
(asterisks) 
and Discrete Leja Points (squares) extracted from the Chebyshev 
lattices on the Lissajous curves, for degree $n=1,2,\dots,30$, compared 
with 
$\mbox{dim}(\mathbb{P}_n^3)=(n+1)(n+2)(n+3)/6$ (upper solid line) 
and $n^2$ (dots).}
\label{leblissa}
\end{figure}
\vskip0.5cm

\vskip0.5cm
\begin{figure}[!ht] \centering
\includegraphics[scale=0.50,clip]{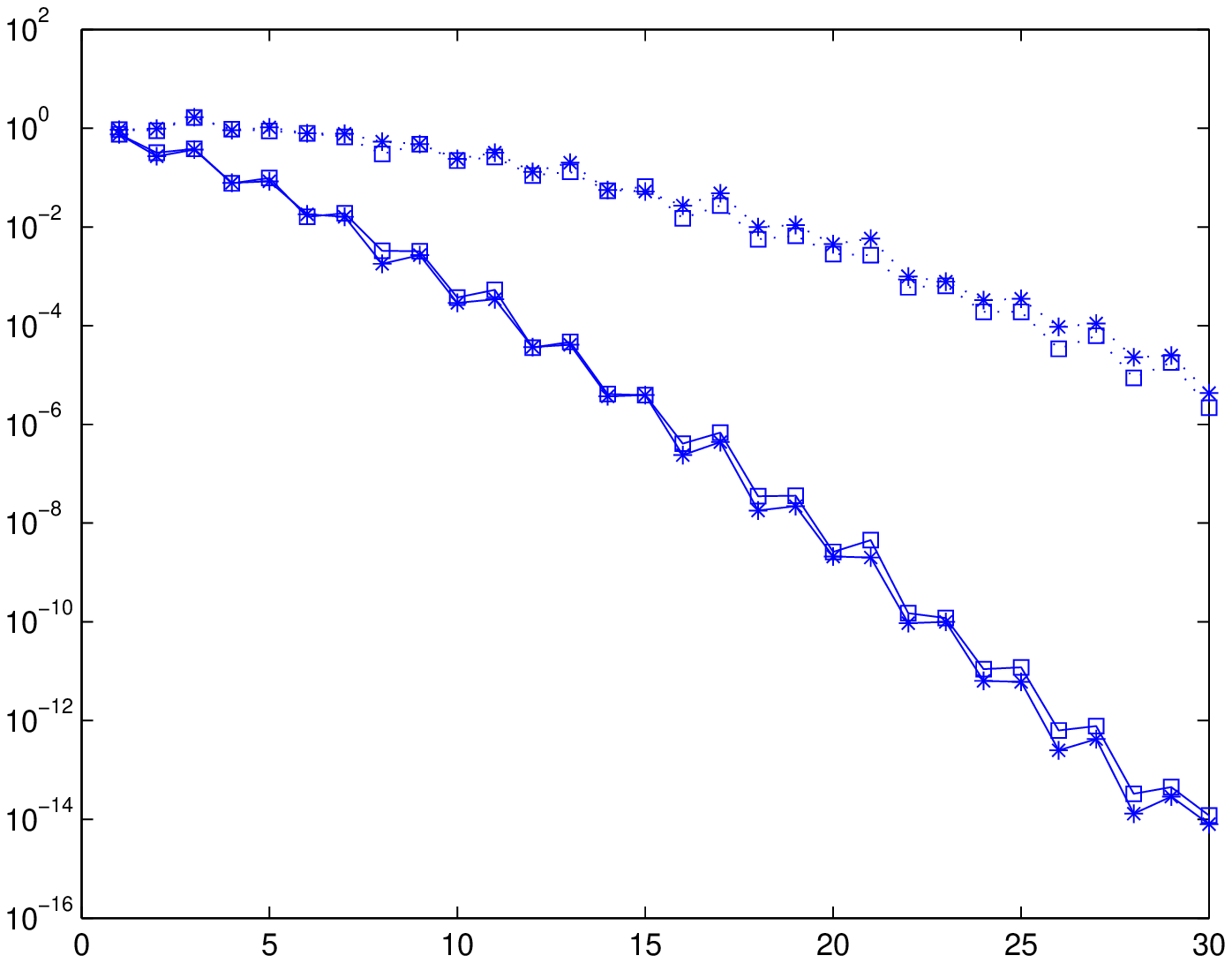}\hfill
\includegraphics[scale=0.50,clip]{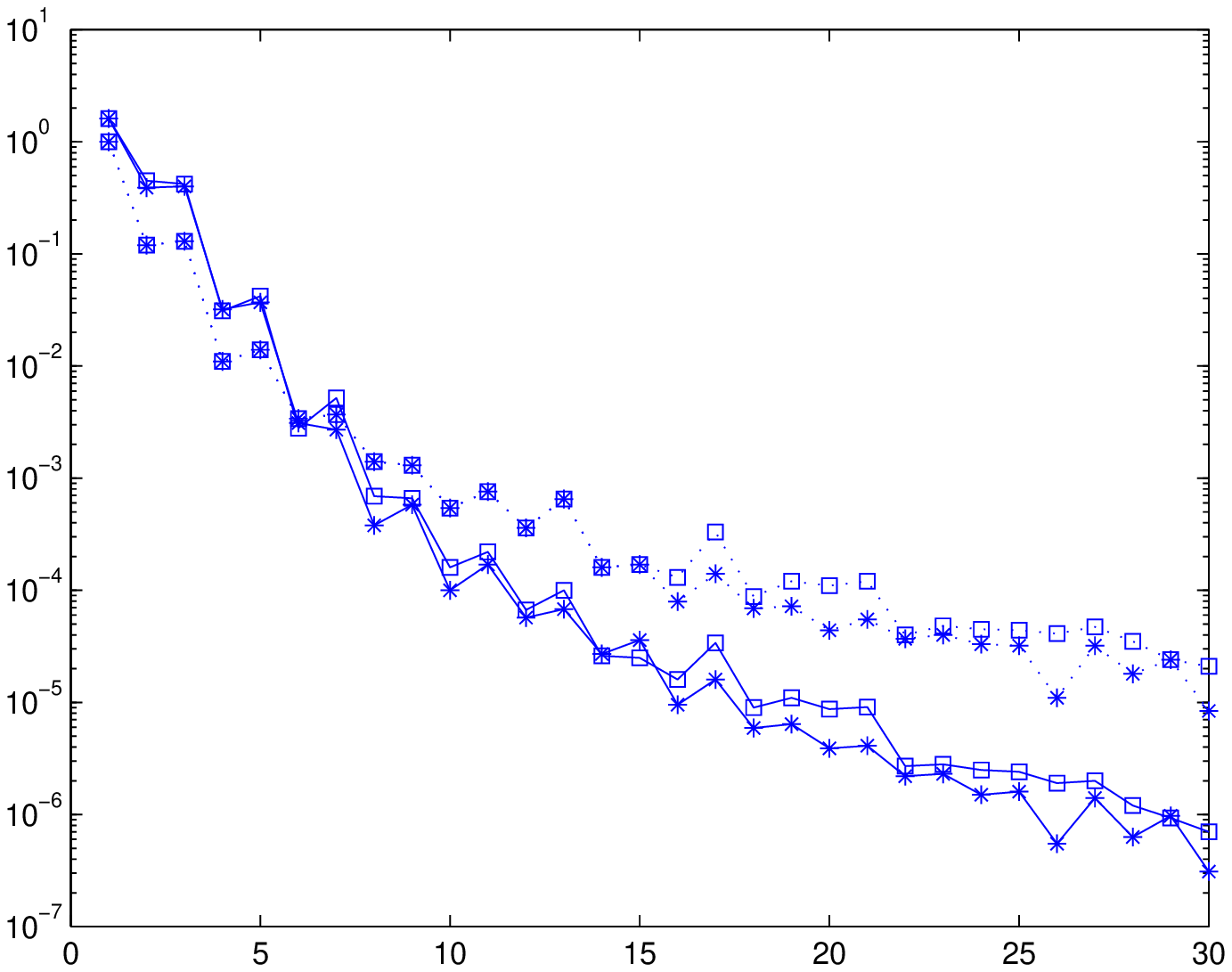}\hfill
\caption{Interpolation errors on Approximate Fekete Points (asterisks) 
and Discrete Leja Points (squares) for the trivariate functions $f_1$ 
(top) 
with $c=1$ (solid line) and $c=5$ (dotted line), and $f_2$ (bottom)   
with $\beta=5$ (solid line) and $\beta=3$ (dotted 
line).}
\label{interperr}
\end{figure}

\section{Appendix}
{\bf Proof of Theorem 1.} We prove the theorem for $n$ even, 
the proof being similar in the odd case. Let be $m=n/2$, $n$ even, 
so that 
$$
(a_n,b_n,c_n)=(3m^2+m,3m^2+2m,3m^2+3m+1)\;. 
$$
We assume that 
$$
i^2+j^2+k^2>0\;.
$$
\noindent
{\bf First case\/}. We show that it is not possible to have 
$$
ia=jb+kc
$$
for $i+j+k\leq 4m$ $(=2n)$. Now, $ia=jb+kc$ becomes 
$i(3m^2+m)+j(3m^2+2m)+k(3m^2+3m)+k$. Since $m$ divides $3m^2+m$, 
$3m^2+2m$ and $3m^2+3m$, we must have that $m$ divides $k$, i.e., 
$k=\alpha m$, $\alpha\geq 0$. Since $k\leq 4m$, $0\leq \alpha\leq 4$. 

Hence we must have   
$$
i(3m^2+m)=j(3m^2+2m)+\alpha m(3m^2+3m+1)
$$
that is, dividing by $m,$
$$
i(3m+1)=j(3m+2)+\alpha(3m^2+3m+1)\;,
$$
which is equivalent to
$$
i((3m+2)-1)=j(3m+2)+\alpha((3m+2)m+(m+1))
$$
and to
$$
(3m+2)(i-j-m\alpha)=i+\alpha(m+1)\;.
$$
The latter implies that 
$$
i+\alpha(m+1)=\beta(3m+2)
$$
for some integer $\beta\geq 0$, i.e.,
$$
i=\beta(3m+2)-\alpha(m+1)
$$
(actually $\beta=i-j-m\alpha$). 

From 
$$
\beta=i-j-m\alpha
$$
we have
$$
j=i-m\alpha-\beta=\beta(3m+2)-\alpha(m+1)-m\alpha-\beta
$$
i.e., 
$$
j=\beta(3m+1)-\alpha(2m+1)
$$
(which must be $\geq 0$). It follows that
$$
i+j+k=\beta(3m+2)-\alpha(m+1)+\beta(3m+1)-\alpha(2m+1)+\alpha m\;,
$$
i.e., 
$$
i+j+k=\beta(6m+3)-\alpha(2m+2)\;.
$$

We now consider two possibilities for $\alpha$:
\begin{itemize}
\item[$1)$] $\alpha=0$. In this case 
$$
i=\beta(3m+2)\;,\;\;j=\beta(3m+1)\;,\;\;k=0
$$
and $i+j+k=\beta(6m+3)$. Now, $\beta\neq 0$ otherwise $i=j=k=0$. 
Hence 
$$
i+j+k\geq 1 (6m+3)>4m
$$
violating the constraint on $i+j+k.$
\item[$2)$] $\alpha\geq 1$ (and $\alpha\leq 4$). In this case $\beta\geq 
1$, for otherwise $i,j<0$. More precisely, since
$$
j=\beta(3m+1)-\alpha(2m+1)=(3\beta-2\alpha)m-\alpha\geq 0
$$ 
we must have $3\beta-2\alpha\geq 1$. Hence
$$
i+j+k=\beta(6m+3)-\alpha(2m+2)=m(6\beta-2\alpha)+3\beta-2\alpha
$$
$$
=m(3\beta-2\alpha+3\beta)+3\beta-2\alpha \geq m(1+3)+1=4m+1>4m
$$
which again violates the constraint on $i+j+k.$
\end{itemize}
\vskip0.5cm
\noindent
{\bf Second case\/}. It is not possible that 
$$
jb=ia+kc
$$ 
for $i+j+k\leq 4m$ $(=2n)$. In this case, $ia=jb+kc$ becomes 
$i(3m^2+m)=j(3m^2+2m)+k(3m^2+3m)+k$. Since $m$ divides $3m^2+m$, $3m^2+2m$ 
and $3m^2+3m$, we must have that $m$ divides $k$, i.e., $k=\alpha m$, 
$\alpha\geq 0$. Since $k\geq 4m$, $0\leq \alpha\leq 4$. 

Hence we must have
$$
j(3m^2+2m)=i(3m^2+m)+\alpha m (3m^2+3m+1)
$$
and dividing by $m$
$$
j(3m+2)=i(3m+1)+\alpha(3m^2+3m+1)
$$ 
which implies that
$$
j(3m+1)+j=i(3m+1)+\alpha(m(3m+1)+2m+1)
$$
and also
$$
j-\alpha(2m+1)=(i-j+\alpha m)(3m+1)\;.
$$

Let $\beta=i-j+\alpha m$ (which a priori could be $\leq 0$) so that 
$$
j-\alpha((2m+1)=\beta(3m+1)
$$
which is equivalent to 
$$
j=\beta(3m+1)+\alpha(2m+1)\;, 
$$
and
$$
i=\beta+j-\alpha m=\beta+(\beta(3m+1)+\alpha (2m+1))-\alpha m\;,
$$
i.e., 
$$
i=\beta(3m+2)+\alpha(m+1)\;.
$$
Hence
\begin{align*}
i+j+k&=\beta(3m+2)+\alpha(m+1)+\beta(3m+1)\\
&\quad+\alpha(2m+1)+\alpha m\\
&=\beta(6m+3)+\alpha(4m+2)\\
&=m(6\beta+4\alpha)+3\beta+2\alpha\\
&=(3\beta+2\alpha)(2m+1)\;.
\end{align*}
For $0<i+j+k\leq 4m$, the only possibility is 
$$
3\beta+2\alpha=1\;.
$$
For $0\leq \alpha\leq 4$, the only integer solution for $\beta$ is 
$$
\alpha=2\;,\;\;\beta=-1\;.
$$
However, in this case, 
$$
i=\beta(3m+2)+\alpha(m+1)=-(3m+2)+2(m+1)=-m<0 
$$
which is not allowed. 
\vskip0.5cm
\noindent
{\bf Third case\/}. It is not possible that 
$$
kc=ia+jb
$$ 
for $i+j+k\leq 4m$ $(=2n)$. In this case, $kc=ia+jb$ becomes
$k(3m^2+3m)+k=i(3m^2+m)+j(3m^2+2m)$. Since $m$ divides $3m^2+m$, 
$3m^2+2m$
and $3m^2+3m$, we must have again that $m$ divides $k$, i.e., $k=\alpha m$,
$\alpha\geq 0$. Since $k\geq 4m$, $0\leq \alpha\leq 4$.

Hence
$$
\alpha m(3m^2+3m+1)=i(3m^2+m)+j(3m^2+2m)\;.
$$
Dividing by $m$ we obtain
$$
\alpha(3m^2+3m+1)=i(3m+1)+j(3m+2)
$$
or equivalently
$$
\alpha(m(3m+2)+m+1)=i(3m+2-1)+j(3m+2)
$$
and 
$$
i+\alpha(m+1)=(3m+2)(-\alpha m+i+j)\;.
$$
Let $\beta=-\alpha m+i+j$. Then
$$
i+\alpha(m+1)=\beta(3m+2)
$$
which implies that
$$
i=\beta(3m+2)-\alpha(m+1)=m(3\beta-\alpha)+(2\beta-\alpha)\;.
$$
Note that $i\geq 0$ implies $\beta\geq 0$ (since $\alpha\geq 0$). 
Further
$$
j=\beta+\alpha m-i=\beta+\alpha m-(\beta(3m+2)-\alpha(m+1))
=\alpha(2m+1)-\beta(3m+1)\;,
$$
i.e., 
$$
j=m(2\alpha-3\beta)+(\alpha-\beta)
$$
and 
$$
i+j+k=\beta(3m+2)-\alpha(m+1)+\alpha(2m+1)-\beta(3m+1)+\alpha m
=\beta+2\alpha m\;.
$$

If $\alpha=0$, then 
$$
i=\beta(3m+2)\;,\;\;j=-\beta(3m+1)\;,\;\;k=0
$$
which is not allowed as $j\ge 0$ (and $\beta\geq 0$).

If $\alpha=3,4$ 
$$
i+j+k=\beta+2\alpha m\geq 6m>4m
$$
which also contradicts the constraints on $i+j+k.$ 

If $\alpha=2$, 
$$
i+j+k=\beta+4m>4m
$$
unless $\beta=0$. However, in this case
$$
i=-2(m+1)<0
$$
and so $\alpha=2$ is not possible. 

The only remaining possibility is $\alpha=1$. In this case
$$
i=\beta(3m+2)-(m+1)\;,\;\;j=(2m+1)-\beta(3m+1)\;,\;\;k=m\;.
$$
But $j\geq 0$ is equivalent to $2m+1\geq \beta(3m+1)$, i.e., 
$$
\beta\leq \frac{2m+1}{3m+1}<1\;,\;\,\mbox{for}\;m\ge1
$$
and so $\beta=0$ (as $\beta$ is an integer). But then 
$$
i=-(m+1)<0
$$
which is not possible. 
\vskip0.5cm
\noindent
{\bf Counterexample\/}. 
Let 
$$
i=2m+1\;,\;\;j=m\;,\;\;k=m\;.
$$
Then $i+j+k=4m+1$ and it is elementary to check that
$ia-jb-kc=0.$
Hence, $4m=2n$ is the maximal value for which the property 
in the statement of Theorem 1 is satisfied.$\;\;\;\square$

\end{document}